\documentclass{article}
\usepackage[utf8]{inputenc}
\usepackage{amssymb,amsmath,amsfonts}
\usepackage[all]{xy}
\usepackage{graphicx}
\usepackage{amsmath}
\usepackage{amsthm}
\usepackage{amsfonts}
\usepackage{amssymb}%
\usepackage{xcolor}
\providecommand{\U}[1]{\protect\rule{.1in}{.1in}}
\newtheorem{theorem}{Theorem}

\newtheorem{cor}[theorem]{Corollary}

\newtheorem{defi}[theorem]{Definition}

\newtheorem{idea memo}[theorem]{Idea Memo}

\newtheorem{notation}[theorem]{Notation}

\newtheorem{remark}[theorem]{Remark}

\mathchardef\mhyphen="2D

\title{Categorical Nonstandard Analysis}
\author{Hayato Saigo\footnote{h\_saigoh@nagahama-i-bio.ac.jp}\\
Nagahama Institute of Bio-Science and Technology\\
Juzo Nohmi\\
Freelance Mathematician}

\begin{document}

\maketitle



\begin{abstract}
In the present paper, we propose a new axiomatic approach to nonstandard analysis and its application to the general theory of spatial structures in terms of category theory. Our framework is based on the idea of internal set theory, while we make use of an endofunctor $\mathcal{U}$ on a topos of sets $S$ together with a natural transformation $\upsilon$, instead of the terms as "standard", "internal" or "external". Moreover, we propose a general notion of a space called $\mathcal{U}$-space, and the category $\mathcal{U}Space$ whose objects are $\mathcal{U}$-spaces and morphisms are functions called $\mathcal{U}$-spatial morphisms. The category $\mathcal{U}Space$, which is shown to be cartesian closed, will give a unified viewpoint toward topological and coarse geometric structure. It will also useful to study  symmetries/asymmetries of the systems with infinite degrees of freedom such as quantum fields.
\end{abstract}

\section{Introduction}\label{sec-Intro}

Nonstandard analysis and category theory are two of the 
great inventions in foundation (or organization) of mathematics . 
Both of them have provided productive viewpoints to organize many kinds of 
topics in mathematics or related fields \cite{Rob, Mac}. 
On the other hand, unification of two theories still seems to be developed, although there are some pioneering works such as \cite{K-M}. 

In the present paper, we propose a new axiomatic framework for nonstandard analysis in terms of category theory.
Our framework is based on the idea of internal set theory \cite{Nel}, while we make use of endofunctor an $\mathcal{U}$ on a topos of sets $\mathcal{S}$ together with a natural transformation $\upsilon$, instead of the terms as "standard", "internal" or "external". 

The triple $(\mathcal{S},\mathcal{U},\upsilon)$ is supposed to satisfy two axioms. The first axiom (``elementarity axiom'') introduced in Section 2 states that the endofunctor $\mathcal{U}$ should preserve all finite limits and finite coproducts. 
Then the endofunctor $\mathcal{U}$ is viewed as some kind of extension of functions preserving all elementary logical properties. In Section 3, we introduce another axiom (``idealization axiom"), which is the translation of ``the principle of idealization'' in internal set theory, and proves the appearance of 
useful entities such as infinitesimals or relations such as ``infinitely close'', in the spirit of Nelson's approach to nonstandard analysis \cite{Nel}. 

Section 4 is devoted to give just a few examples of applications on topology (on metric spaces, for simplicity). Although the characterizations of continuous maps or uniform continuous maps in terms of nonstandard analysis are well known, but we prove them from our framework for the reader's convenience. In section 5, we characterize the notion of bornologous map, which is a fundamental notion in coarse geometry \cite{Roe}. 

In section 6, we introduce the notion of $\mathcal{U}$-space and $\mathcal{U}$-morphism, which are the generalization of examples in the previous two sections. We introduce the category $\mathcal{U}Space$ consisting of $\mathcal{U}$-spaces and $\mathcal{U}$-morphisms, which is shown to be cartesian closed. It will give a unified viewpoint toward topological and coarse geometric structure, and will be useful to study symmetries/asymmetries of the systems with infinite degrees of freedom such as quantum fields.

\section{Elementarity Axiom}

Let $\mathcal{S}$ be a topos of sets, i.e., an elementary topos with natural number object satisfying well-pointedness and the axiom of choice (See \cite{Mac}, which is based on the idea in \cite{ETCS}).
We will make use of an endofunctor $\mathcal{U}:\mathcal{S}\longrightarrow \mathcal{S}$ with a natural transformation $\upsilon: Id_{\mathcal{S}}\longrightarrow\mathcal{U}$ satisfying two axioms, ``elementarity axiom'' and ``idealization axiom''.

\begin{quote}
\textbf{Elementarity Axiom}:
$\mathcal{U}$ preserves all finite limits and finite coproducts.
\end{quote}

\begin{remark}
$\mathcal{U}$ does not necessarily preserve power sets. This is the reason for the name of ``elementarity''.
\end{remark}


It is easy to see that ``elementarity axiom'' implies preservation of many basic notions such as elements, subsets, finite cardinals (especially the subobject classifier 2) and propositional calculi. Moreover, 
the following theorem holds.

\begin{theorem}
$\mathcal{U}$ is faithful.
\begin{proof}
Because it preserves diagonal morphisms and complements.
\end{proof}
\end{theorem}

\begin{theorem}\label{thm:iota}
For any element $x:1\longrightarrow X$, $\upsilon_X (x)=\mathcal{U}(x) \circ \upsilon_1$.
\begin{proof}
By naturality of $\upsilon$.
\end{proof}
\end{theorem}
\begin{cor}
All components of $\upsilon$ are monic. 
\end{cor}
From the discussion above, a set $X$ in $\mathcal{S}$ is to be considered as a canonical subset of $\mathcal{U}(X)$
through $\upsilon_X:X\longrightarrow \mathcal{U}(X)$. Hence, $\mathcal{U}(f): \mathcal{U}(X_0)\longrightarrow \mathcal{U}(X_1)$ 
can be considered as ``the function induced from $f:X_0\longrightarrow X_1$ through $\upsilon$.'' 


\begin{defi}
Let $A,B$ be objects in $\mathcal{S}$. 
The function $ev_{A,B}: A\times B^A \longrightarrow B$ satisfying
\[
ev_{A,B}(a,f)=f(a)
\]
for all $(a,f)\in A\times B^A$ is called the evaluation (for $A,B$) . The lambda conversion $\hat{g}:Z\longrightarrow Y^X$ of $g:X\times Z\longrightarrow Y$ is the function satisfying
\[
g=ev_{X,Y}\circ (1_X\times \hat{g}),
\] 
where $1_X\times \hat{g}$ denotes the function satisfying
\[
(1_X\times \hat{g})(x,z)=(x,\hat{g}(z)).
\]
We define a family of functions 
$\kappa_{A,B}:\mathcal{U}(B^A)\longrightarrow \mathcal{U}(B)^{\mathcal{U}(A)}$ in $\mathcal{S}$ by 
the lambda conversion of $\mathcal{U}(ev_{A,B}):\mathcal{U}(A\times B^A)\cong\mathcal{U}(A)\times \mathcal{U}(B^A)\longrightarrow \mathcal{U}(B)$. 

\end{defi}


The theorem below means that $\kappa_{A,B}\circ \upsilon_{B^A}$ represents ``inducing $\mathcal{U}(f)$ from $f$ through $\upsilon$'' 
in terms of exponentials.
\begin{theorem}
Let $f:A\longrightarrow B$ be any function in $\mathcal{S}$. Then
\[
\kappa_{A,B}\circ \upsilon_{B^A}(\widehat{f})=\widehat{\mathcal{U}(f)}.
\]
(Here ~$\:\widehat{}$~ denotes the lambda conversion operation.)  
\begin{proof}
Take the (inverse) lambda conversion 
of the left hand side of the equality to be proved. 
It is 
$ev_{\mathcal{U}(A),\mathcal{U}(B)} \circ
(id_{\mathcal{U}(A)} \times \kappa_{A,B}) \circ
(id_{\mathcal{U}(A)} \times \upsilon _{B^A}) \circ
(id_{\mathcal{U}(A)} \times \widehat{f})$.
By naturality of $\upsilon$ and functorial properties of $\mathcal{U}$, it is calculated as follows:
\begin{align*}
\xymatrix@C=1.7cm{
& \mathcal{U} (A) \times B^A \ar[rd]^-{id_{\mathcal{U}(A)} \times \upsilon _{B^A}} \\
\mathcal{U}(A) \times 1 \ar[r]_-{\sim}
\ar[ru]^-{id_{\mathcal{U}(A)} \times \widehat{f}}
\ar[d]_-{\rotatebox{90}{$\sim$}} &
\mathcal{U}(A) \times \mathcal{U}(1) \ar[r]_{id_{\mathcal{U}(A)} \times \mathcal{U}(\widehat{f})}
\ar[d]_-{\rotatebox{90}{$\sim$}} &
\mathcal{U}(A) \times \mathcal{U}(B^A) \ar[rd]^-{id_{\mathcal{U}(A)} \times \kappa_{A,B}}
\ar[d]_-{\rotatebox{90}{$\sim$}} \\
\mathcal{U}(A) \ar[r]_-{\sim} \ar[drr]_{\mathcal{U}(f)} &
\mathcal{U}(A \times 1) \ar[r]_-{\mathcal{U}(id_{A} \times \widehat{f})} &
\mathcal{U}(A \times B^A) \ar[d]_{\mathcal{U}(ev_{A,B})} &
\mathcal{U}(A) \times \mathcal{U}(B)^{\mathcal{U}(A)}
\ar[ld]^-{~~~~~~ev_{\mathcal{U}(A),\mathcal{U}(B)}} \\
&&\mathcal{U}(B)
}
\end{align*}
\end{proof}
\end{theorem}

\begin{cor}
$\kappa_{A,B}\circ \upsilon_{B^A}$ is monic. 
\end{cor}

\begin{notation}
\textbf{From now on, we omit $\upsilon$ and $\kappa$. $\mathcal{U}(f):\mathcal{U}(X)\longrightarrow \mathcal{U}(Y)$ will be often identified with $f:X\longrightarrow Y$ and denoted simply 
as $f$ instead of $\mathcal{U}(f)$.}

\end{notation}

\begin{theorem}

Let $P:X\longrightarrow 2$ be any proposition (function in $\mathcal{S}$). Then
\[
\forall_{x\in X} P(x) \Longleftrightarrow \forall_{x\in \mathcal{U}(X)} P(x).
\]
\begin{proof}
$P:X\longrightarrow 2$ factors thorough $``true'': 1\longrightarrow 2$ if and only if $P:\mathcal{U}(X)\longrightarrow 2$ factors thorough $``true'': 1\longrightarrow 2$.
\end{proof} 
\end{theorem}

Dually, we obtain the following:

\begin{theorem}

Let $P:X\longrightarrow 2$ be any proposition (function in $\mathcal{S}$). Then
\[
\exists_{x\in X} P(x) \Longleftrightarrow \exists_{x\in \mathcal{U}(X)} P(x).
\]

\end{theorem}

The two theorems above are considered as the simplest versions of ``transfer principle''. 
To treat with free variables and quantification, the theorem below is important. (The author thanks Professor Anders Kock for indicating this crucial point.)

\begin{theorem}
$\mathcal{U}$ preserves images.

\begin{proof}
As $\mathcal{U}$ preserves all finite limits, it preserves monics. On the other hand, it also preserves epic since every functor preserves split epics and every epic in $\mathcal{S}$ is split epic (axiom of choice). 
Hence the image, which is nothing but the epi-mono factorization, is preserved. 
\end{proof}
\end{theorem}


\section{Idealization Axiom}
From our viewpoint, nonstandard Analysis is nothing but a method of using an endofunctor which satisfies ``elementarity axiom'' and the following 
``idealization axiom''. The name is after ``the principle of idealization'' in Nelson's internal set theory \cite{Nel}. Most of the basic idea 
in this section is much in common with \cite{Nel}, although the functorial approach is not taken in internal set theory.

\begin{remark}
Internal set theory (IST) is a syntactical approach to nonstandard analysis consisting of the ``principle of Idealization (I)'' and the two more 
basic principles called  ``principle of Standardization (S)'' and ``Transfer principle (T)''.
 In our framework, the role of (S) is played by the axiom of choice for $\mathcal{S}$ and 
(T) corresponds to the contents of section 2.     
\end{remark}

\begin{notation}
For any set $X$, $\tilde{X}$ denotes the set of all finite subsets of $X$. 
\end{notation}

\begin{quote}

\textbf{Idealization Axiom}:
Let $P$ be an element of $\mathcal{U}(2^{X\times Y})$. Then
\[
\forall_{x'{\in \tilde{X}}}\exists_{y\in \mathcal{U}(Y)}\forall_{x\in x'} P(x,y)\Longleftrightarrow \exists_{y\in \mathcal{U}(Y)} \forall_{x\in X}P(x,y) .
\]
   
\end{quote}
Or dually,
\begin{quote}

\textbf{Idealization Axiom, dual form}:
Let $P$ be an element of $\mathcal{U}(2^{X\times Y})$
Then
\[
\exists_{x'{\in \tilde{X}}}\forall_{y\in \mathcal{U}(Y)}\exists_{x\in x'} P(x,y)\Longleftrightarrow \forall_{y\in \mathcal{U}(Y)} \exists_{x\in X}P(x,y) .
\]

\end{quote}




When $X$ is a directed set with an order $\leq $ and $P\in \mathcal{U}(2^{X\times Y})$ satisfies the ``filter condition'', i.e.,
\[
\forall_{x_0 \in X} ( P(x_0,y) \Longrightarrow \forall_{x \in X}(x\leq x_0\Longrightarrow P(x,y)) ),
\]
or dually, the ``cofilter condition'', i.e.,
\[
\forall_{x_0 \in X} ( P(x_0,y) \Longrightarrow \forall_{x \in X} (x_0\leq x\Longrightarrow P(x,y)) ),
\]
then idealization axiom is simplified as ``commutation principle'':

\begin{theorem}[Commutation Principle]
If $P\in \mathcal{U}(2^{X\times Y})$ satisfies the ``filter condition'' and ``cofilter condition'' above, respectively, then
\[
\forall_{x\in X}\exists_{y\in \mathcal{U}(Y)}P(x,y)\Longleftrightarrow \exists_{y\in \mathcal{U}(Y)} \forall_{x\in X}P(x,y)
\]
and
\[
\exists_{x\in X}\forall_{y\in \mathcal{U}(Y)} P(x,y)\Longleftrightarrow \forall_{y\in \mathcal{U}(Y)} \exists_{x\in X}P(x,y)
\]
holds, respectively.
\end{theorem}

By the principle above, we can easily prove the existence of ``unlimited numbers'' in $\mathcal{U}(\mathbb{N})$, where all arithmetic operation and order structure on $\mathbb{N}$ is naturally extended.

\begin{theorem}[Existence of ``unlimited numbers'']

There exists some $\omega \in \mathcal{U}(\mathbb{N})$ such that $n \leq \omega$ for any $n\in \mathbb{N}$. 

\begin{proof}
Because it is obvious that for any $n \in \mathbb{N}$ there exists some $\omega \in \mathbb{N}\subset \mathcal{U}(\mathbb{N})$ such that $n< \omega $.
\end{proof}
\end{theorem}

As in $\mathcal{S}$ we can construct rational numbers and the completion of them as usual, we have the object $\mathbb{R}$, the set of real numbers. Then we obtain the following:

\begin{cor}
``Infinitesimals'' do exist in $\mathcal{U}(\mathbb{R})$. That is, there exists some $r \in \mathcal{U}(\mathbb{R})$ such that $|r|<R$ for any positive $R \in \mathbb{R}$. 
\end{cor}

\section{Topological Structure: Continuous Map and Uniform Continuous Map}

We will take an example of basic applications of nonstandard analysis within our framework, i.e., characterization of continuity and uniform continuity in terms of a relation $\approx$ (``infinitely close'') 
on $\mathcal{U}(X)$, which is based on essentially the same arguments well-known in nonstandard analysis, especially, internal set theory\cite{Nel}. For simplicity, we will discuss only for metric spaces here. (For more general topological spaces, we can define $\approx$ in terms of the system of open sets. See \cite{Nel} for example.)

\begin{defi}[Infinitely close]
Let $(X, d)$ be a metric space. We call the relation $\approx$ on $\mathcal{U}(X)$ defined below as ``infinitely close'':
\[
x\approx x' \Longleftrightarrow \forall_{\epsilon\in \mathbb{R}} d(x,x')<\epsilon.
\]

\end{defi}
That is, $d(x,x')$ is infinitesimal. It is easy to see that $\approx$ is an equivalence relation on $\mathcal{U}(X)$. 

\begin{theorem}[Characterization of continuity]
Let $(X_0,d_0), (X_1, d_1)$ be metric spaces and $\approx_0, \approx_1$ be infinitely close relations on them, respectively. A map $f:X_0\longrightarrow X_1$ is continuous if and only if
\[
\forall_{x\in X_0}\forall_{x'\in\mathcal{U}(X_0)} \: (\:  x\approx_0 x' \Rightarrow f(x)\approx_1 f(x')\: )
\]
holds.
\begin{proof} We can translate the condition for $f$ by using usual logic, ``commutation principle'' and ``transfer principle'' as follows:
\begin{eqnarray*}
 &\forall_{x\in X_0}\forall_{x'\in\mathcal{U}(X_0)} \: (\:  x\approx_0 x' \Rightarrow f(x)\approx_1 f(x')\: )& \\
 \Longleftrightarrow & \forall_{x\in X_0}\forall_{x'\in\mathcal{U}(X_0)}  \:( \: \forall_{\delta \in \mathbb{R}} \: d_0(x,x')<\delta \Rightarrow \forall_{\epsilon \in \mathbb{R}}\: d_1(f(x), f(x'))<\epsilon \:)& \\
 \Longleftrightarrow & \forall_{x\in X_0}\forall_{x'\in\mathcal{U}(X_0)} \forall_{\epsilon \in \mathbb{R}} \exists_{\delta \in \mathbb{R}}\:  (\: d_0(x,x')<\delta \Rightarrow d_1(f(x), f(x'))<\epsilon \:) & \\
 \Longleftrightarrow & \forall_{x\in X_0} \forall_{\epsilon \in \mathbb{R}} \forall_{x'\in\mathcal{U}(X_0)}  \exists_{\delta \in \mathbb{R}}\:  (\: d_0(x,x')<\delta \Rightarrow d_1(f(x), f(x'))<\epsilon \:) & \\
 \Longleftrightarrow & \forall_{x\in X_0}\forall_{\epsilon \in \mathbb{R}} \exists_{\delta \in \mathbb{R}} \forall_{x'\in\mathcal{U}(X_0)}\:   (\: d_0(x,x')<\delta \Rightarrow d_1(f(x), f(x'))<\epsilon \:) & \\
 \Longleftrightarrow & \forall_{x\in X_0}\forall_{\epsilon \in \mathbb{R}} \exists_{\delta \in \mathbb{R}}  \forall_{x'\in X_0}\:   (\: d_0(x,x')<\delta \Rightarrow d_1(f(x), f(x'))<\epsilon \:). & 
\end{eqnarray*}
\end{proof}
\end{theorem}

\begin{theorem}[Characterization of uniform continuity]
Let $(X_0,d_0), (X_1, d_1)$ be metric spaces and $\approx_0, \approx_1$ be infinitely close relations on them, respectively. A map $f:X_0\longrightarrow X_1$ is uniformly continuous if and only if
\[
\forall_{x\in\mathcal{U}(X_0)}\forall_{x'\in\mathcal{U}(X_0)} \: (\:  x\approx_0 x' \Rightarrow f(x)\approx_1 f(x')\: )
\]
holds.
\begin{proof} We can translate the condition for $f$ by using usual logic, ``commutation principle'' and ``transfer Principle'' as follows:
\begin{eqnarray*}
 &\forall_{x\in\mathcal{U}(X_0)}\forall_{x'\in\mathcal{U}(X_0)} \: (\:  x\approx_0 x' \Rightarrow f(x)\approx_1 f(x')\: )& \\
 \Longleftrightarrow & \forall_{x\in\mathcal{U}(X_0)}\forall_{x'\in\mathcal{U}(X_0)}  \:( \: \forall_{\delta \in \mathbb{R}} \: d_0(x,x')<\delta \Rightarrow \forall_{\epsilon \in \mathbb{R}}\: d_1(f(x), f(x'))<\epsilon \:)& \\
 \Longleftrightarrow & \forall_{x\in\mathcal{U}(X_0)}\forall_{x'\in\mathcal{U}(X_0)} \forall_{\epsilon \in \mathbb{R}} \exists_{\delta \in \mathbb{R}}\:  (\: d_0(x,x')<\delta \Rightarrow d_1(f(x), f(x'))<\epsilon \:) & \\
 \Longleftrightarrow & \forall_{\epsilon \in \mathbb{R}} \forall_{x\in\mathcal{U}(X_0)}\forall_{x'\in\mathcal{U}(X_0)}  \exists_{\delta \in \mathbb{R}}\:  (\: d_0(x,x')<\delta \Rightarrow d_1(f(x), f(x'))<\epsilon \:) & \\
 \Longleftrightarrow & \forall_{\epsilon \in \mathbb{R}} \exists_{\delta \in \mathbb{R}} \forall_{x\in\mathcal{U}(X_0)}\forall_{x'\in\mathcal{U}(X_0)}\:   (\: d_0(x,x')<\delta \Rightarrow d_1(f(x), f(x'))<\epsilon \:) & \\
 \Longleftrightarrow & \forall_{\epsilon \in \mathbb{R}} \exists_{\delta \in \mathbb{R}} \forall_{x\in X_0} \forall_{x'\in X_0}\:   (\: d_0(x,x')<\delta \Rightarrow d_1(f(x), f(x'))<\epsilon \:). & 
\end{eqnarray*}
\end{proof}
\end{theorem}

 As we have seen, a morphism between metric spaces is characterized as ``a morphism with respect to $\approx$''. This suggests the possibility for considering 
 other kinds of ``equivalence relations on (some subset of) $\mathcal{U}(X)$'' as generalized spatial structures on $X$. 
  In the next section, we will take one example related to  large scale geometric structure.

\section{Coarse Structure: Bornologous Map}


Let us consider another kind of equivalence relation $\sim$ (``finitely remote'') defined below. For simplicity, we will discuss only for metric spaces here.

\begin{defi}[Finitely remote]
Let $(X, d)$ be a metric space. We call the relation $\sim$ on $\mathcal{U}(X)$ defined below as ``finitely remote'':
\[
x\sim x' \Longleftrightarrow \exists_{R \in \mathbb{R}} d(x,x')<R.
\]
\end{defi}

Note that we use $\exists$ instead of $\forall$, in contrast to ``infinitely close''. This kind of dual viewpoint will be proved to be useful in the geometric study of large scale structure such as coarse geometry\cite{Roe}.  

In fact, we can prove that ``bornologous map'', a central notion of a morphism for coarse geometry, can be characterized as ``a morphism with respect to $\sim$'', just like (uniform) continuity can be viewed as ``a morphism with respect to $\approx$''.

\begin{defi}[Bornologous map]
Let $(X_0,d_0)$ and $(X_1, d_1)$ be metric spaces. A map $f:X_0\longrightarrow X_1$ is called a bornologous map when 
\[
\forall_{R \in \mathbb{R}} \exists_{S \in \mathbb{R}}\: \forall_{x\in X_0} \forall_{x'\in X_0}   (\: d_0(x,x')<R \Rightarrow d_1(f(x), f(x'))<S \:)
\]
holds.
\end{defi}

\begin{theorem}[Characterization of bornologous map]
Let $(X_0,d_0), (X_1, d_1)$ be metric spaces and $\sim_0, \sim_1$ be finitely remote relations on them, respectively. A map $f:X_0\longrightarrow X_1$ is bornologous if and only if
\[
\forall_{x\in\mathcal{U}(X_0)}\forall_{x'\in\mathcal{U}(X_0)} \: (\:  x\sim_0 x' \Rightarrow f(x)\sim_1 f(x')\: )
\]
holds.
\begin{proof} We can translate the condition for $f$ by using usual logic, ``commutation principle'' and ``transfer principle'' as follows:
\begin{eqnarray*}
 &\forall_{x\in\mathcal{U}(X_0)}\forall_{x'\in\mathcal{U}(X_0)} \: (\:  x\sim_0 x' \Rightarrow f(x)\sim_1 f(x')\: )& \\
 \Longleftrightarrow & \forall_{x\in\mathcal{U}(X_0)}\forall_{x'\in\mathcal{U}(X_0)}  \:( \: \exists_{R \in \mathbb{R}} \: d_0(x,x')<R \Rightarrow \exists_{S \in \mathbb{R}}\: d_1(f(x), f(x'))<S \:)& \\
 \Longleftrightarrow & \forall_{x\in\mathcal{U}(X_0)}\forall_{x'\in\mathcal{U}(X_0)} \forall_{R \in \mathbb{R}} \exists_{S \in \mathbb{R}}\:  (\: d_0(x,x')<R \Rightarrow d_1(f(x), f(x'))<S \:) & \\
 \Longleftrightarrow & \forall_{R \in \mathbb{R}} \forall_{x\in\mathcal{U}(X_0)}\forall_{x'\in\mathcal{U}(X_0)}  \exists_{S \in \mathbb{R}}\:  (\: d_0(x,x')<R \Rightarrow d_1(f(x), f(x'))<S \:) & \\
 \Longleftrightarrow & \forall_{R \in \mathbb{R}} \exists_{S \in \mathbb{R}} \forall_{x\in\mathcal{U}(X_0)}\forall_{x'\in\mathcal{U}(X_0)}\:   (\: d_0(x,x')<R \Rightarrow d_1(f(x), f(x'))<S \:) & \\
 \Longleftrightarrow & \forall_{R \in \mathbb{R}} \exists_{S \in \mathbb{R}} \forall_{x\in X_0} \forall_{x'\in X_0}\:   (\: d_0(x,x')<R \Rightarrow d_1(f(x), f(x'))<S \:). & 
\end{eqnarray*}
\end{proof}
\end{theorem}


\section{The notion of $\mathcal{U}$-space and the Category $\mathcal{U} Space$}

Based on the characterizations of topological and coarse geometrical structure, we introduce the notion of $\mathcal{U} $-space.

\begin{defi} [$\mathcal{U} $-space]
A $\mathcal{U} $-space is a triple $(X, K, \rightsquigarrow )$ consisting of a set $X$, a subset $K$ of $\mathcal{U}(X)$ which includes $X$ as a subset and a preorder $\rightsquigarrow $ defined on $K$.
\end{defi}

When the preorder $\rightsquigarrow $ is an equivalence relation, i.e., a preorder satisfying symmetry, we call the $\mathcal{U} $-space symmetric. A symmetric $\mathcal{U} $-space $(X, K, \rightsquigarrow )$ is called uniform if $K=\mathcal{U}(X)$. The "infinitely close" relation and the "finitely remote" relation provide the simplest examples of uniform $\mathcal{U}$-space structure.

    

Actually, any topological space $X$ with the set of open sets $T$ can be viewed as $\mathcal{U}$-space $(X,\mathcal{U}(X),\rightharpoonup)$ where  $x\rightharpoonup x'$ denotes the preorder "$\forall O\in T\:\: x\in \mathcal{U}(O)\Longrightarrow x'\in \mathcal{U}(O)$". If $(X,T)$ is a Hausdorff space, we can construct the symmetric $\mathcal{U}$-space $(X,K,\rightsquigarrow)$, where $K$ denotes 
\[
K=\{x\in \mathcal{U}(X)| \exists x_0\in X \: x\rightharpoonup x'\}
\]
and $x\rightsquigarrow x'$ is defined as the relation "$\exists x_0\in X \:\: x_0 \rightharpoonup x \& x_0 \rightharpoonup x'$." The transitivity of $\rightharpoonup$ follows from the fact that if $(X,T)$ is Hausdorff, $x_0\rightharpoonup x$ and $x'_0\rightharpoonup x$ imply $x_0=x'_0$ for all $x_0,x'_0\in X$. Actually, the preorder $\rightsquigarrow$ becomes an equivalence relation.  

The concept of $\mathcal{U}$-space will provide a general framework to unify various spatial structure such as topological structure and coarse structure. 
The notion of morphism between $\mathcal{U}$-spaces is defined as follows:

\begin{defi}[$\mathcal{U}$-spatial morphism] Let $(X_0,K_0,\rightsquigarrow_0)$ and $(X_1,K_1, \rightsquigarrow_1)$ be $\mathcal{U}$-spaces. A function $f:X_0\rightarrow X_1$ is called a $\mathcal{U}$-spatial morphism from 
$(X_0,K_0,\rightsquigarrow_0)$ to $(X_1,K_1, \rightsquigarrow_1)$ when $f(K_0)\subset K_1$ and 
\[
x\rightsquigarrow_0 x' \Longrightarrow f(x)\rightsquigarrow_1 f(x')
\] 
holds for any $x,x'\in K_0$.
\end{defi}

The uniform continuous maps and bornologous maps between metric spaces are nothing but $\mathcal{U}$-spatial morphisms between corresponding uniform $\mathcal{U}$-spaces. The notion of continuous maps between Hausdorff spaces can be characterized as $\mathcal{U}$-spatial morphisms between the corresponding symmetric $\mathcal{U}$-spaces. 


\begin{defi}[Category $\mathcal{U} Space$]
The category $\mathcal{U} Space$ is a category whose objects are $\mathcal{U}$-spaces and whose morphisms are $\mathcal{U}$-spatial morphisms.
\end{defi}

\begin{defi}Let $(X_0,K_0,\rightsquigarrow_0)$ and $(X_1,K_1, \rightsquigarrow_1)$ be $\mathcal{U}$-spaces. The $\mathcal{U}$-space $(X_0\times X_1, K_0\times K_1, \rightsquigarrow)$, where the preorder $\rightsquigarrow$ is defined as
\[
(x_0,x_1)\rightsquigarrow (x'_0,x'_1) \Longleftrightarrow x_0\rightsquigarrow_0 x'_0\: \&\:  x_1\rightsquigarrow_1 x'_1,
\]
is called the product $\mathcal{U}$-space of $(X_0,K_0,\rightsquigarrow_0)$ and $(X_1,K_1, \rightsquigarrow_1)$.
\end{defi}

\begin{theorem}
The projections become $\mathcal{U}$-spatial morphisms. The diagram consisting of two $\mathcal{U}$-spaces, the product space of them and projections become a product in $\mathcal{U} Space$. 
\begin{proof}
Easy.
\end{proof}
\end{theorem}

\begin{defi}[Exponential $\mathcal{U}$-space]

Let $(X_0,K_0,\rightsquigarrow_0)$ and $(X_1, K_1,\rightsquigarrow_1)$ be $\mathcal{U}$-spaces. We denote the set of all $\mathcal{U}$-spatial morphisms from $X_0$ to $X_1$ as $[X_1^{X_0}]$, which is the subset of $X_1^{X_0}$. 
The restriction of $ev_{X_0,X_1}:X_0\times X_1^{X_0}\longrightarrow X_1$ onto $X_0\times [X_1^{X_0}]\longrightarrow X_1$ is denoted as $[ev_{X_0,X_1}]$. The $\mathcal{U}$-space $([X_1^{X_0}],K,\rightsquigarrow)$, where $K$ is defined as the subset of $\mathcal{U}([X_1^{X_0}])$,
\[
K=\{f
| 
\forall x\in K_0 [ev_{X_0,X_1}](x,f)\in K_1 \; \text{and}\; x\rightsquigarrow_0  x'\Longrightarrow[ev_{X_0,X_1}](x,f)\rightsquigarrow_1 [ev_{X_0,X_1}](x,f')
\}
\]
and $\rightsquigarrow$ is defined as
\[
f\rightsquigarrow f' \Longleftrightarrow \forall x \in K_0\;\;[ev_{X_0,X_1}](x,f)\rightsquigarrow_1 [ev_{X_0,X_1}](x,f'), 
\]
is called the exponential $\mathcal{U}$-space from $(X_0,K_0,\rightsquigarrow_0)$ to $(X_1, K_1,\rightsquigarrow_1)$.
\end{defi}

\begin{theorem}
Let $(X_0,K_0,\rightsquigarrow_0)$ and $(X_1,K_1,\rightsquigarrow_1)$ be $\mathcal{U}$-spaces and $([X_1^{X_0}],K,\rightsquigarrow)$ be the exponential $\mathcal{U}$-space from $(X_0,K_0,\rightsquigarrow_0)$ to $(X_1,K_1,\rightsquigarrow_1)$.
The morphism $[ev_{X_0,X_1}]:X_0\times [X_1^{X_0}]\longrightarrow X_1$, the restriction of $ev_{X_0,X_1}$, is a $\mathcal{U}$-spatial morphism.
Moreover, it becomes an evaluation in $\mathcal{U} Space$ and $[X_1^{X_0}]$ is an exponential in $\mathcal{U} Space$.
\begin{proof}
First we prove that $[ev_{X_0,X_1}]$ is a $\mathcal{U}$-spatial morphism:
For any $(x,f)\in K_0\times K$, $[ev_{X_0,X_1}](x,f)$ is in $K_1$ since $x\in K_0$ and $f\in K$. Suppose that $(x,f),(x',f')\in K_0\times K$ and  $(x,f)\rightsquigarrow (x',f')$, that is, $x,x'\in K_0$, $f,f'\in K$, $x\rightsquigarrow x'$ and $f\rightsquigarrow f'$. Then we have 
\[
[ev_{X_0,X_1}](x,f)\rightsquigarrow [ev_{X_0,X_1}](x,f') 
\] 
since $f\rightsquigarrow f'$. We also have
\[
[ev_{X_0,X_1}](x,f')\rightsquigarrow [ev_{X_0,X_1}](x',f')
\]
since $f' \in K$. 
Hence, $[ev_{X_0,X_1}](x,f)\rightsquigarrow [ev_{X_0,X_1}](x',f')$. 

Next we prove that $[ev_{X_0,X_1}]:X_0\times [X_1^{X_0}]\longrightarrow X_1$ becomes an evaluation in $\mathcal{U} Space$ and $[X_1^{X_0}]$ is an exponential in $\mathcal{U} Space$:
Let $(X_2,K_2,\rightsquigarrow_2)$ be any $\mathcal{U}$-space and $f:X_0\times X_2\longrightarrow X_1$ be any $\mathcal{U}$-spatial morphism. 
Consider the lambda conversion $\hat{f}:X_2\longrightarrow X_1^{X_0}$. 
By assumption that $f$ is $\mathcal{U}$-spatial,
\[
(x,c),(x',c') \in K_0\times K_2 \;\& \; (x,c)\rightsquigarrow (x',c') \Longrightarrow f(x,c),f(x',c') \in K_1 \;\&\; f(x,c)\rightsquigarrow_1 f(x',c')
\]
holds,
where $\rightsquigarrow$ denote the preorder on $K_0\times K_2$.
It is equivalent to the statement that
$c,c'\in K_2$ and $c\rightsquigarrow_2 c'$ imply 
\[
x,x'\in K_0 \& x\rightsquigarrow_0 x'\Longrightarrow \hat{f}(c)(x),\hat{f}(c')(x')\in K_1 \& \hat{f}(c)(x)\rightsquigarrow_1 \hat{f}(c')(x').
\]
Applying the implication above for the case $c=c'\in X_2$, we have $\hat{f}(c)\in [X_1^{X_0}]$. Hence, we can replace $\hat{f}:X_2\longrightarrow X_1^{X_0}$ with $[\hat{f}]:X_2\longrightarrow [X_1^{X_0}]$ by restricting the codomain to the image of $\hat{f}$.
Moreover, we can also prove that $[\hat{f}]$ is $\mathcal{U}$-spatial from the implication: By the implication above we have $[\hat{f}](c),[\hat{f}](c')\in K$ and $[\hat{f}](c) \rightsquigarrow [\hat{f}](c')$ when $c,c'\in K_2$ and $c\rightsquigarrow_2  c'$. This means that $[\hat{f}]$ is $\mathcal{U}$-spatial.

It is easy to show that this $[\hat{f}]$ is the unique $\mathcal{U}$-spatial morphism from $X_2$ to $[X_1^{X_0}]$ satisfying $[ev_{X_0,X1}]\circ (1_{X_0}\times [\hat{f}])=f$. That completes the proof.
\end{proof}
\end{theorem}

Combining the two theorems above, we have:

\begin{theorem}
The category $\mathcal{U} Space$ is a cartesian closed category.
\end{theorem}

\section*{Acknowledgements}
The authors would like to express his sincere thanks to Prof. Anders Kock and Prof. Edward Nelson for their encouragements. They are grateful to Prof. Hiroshi Ando, Prof. Izumi Ojima, Dr. Kazuya Okamura, Ms. Misa Saigo, Prof. Hiroki Sako and Prof. Ryokichi Tanaka for fruitful discussions and comments.
This work  was partially supported by Research Origin for Dressed
Photon, JSPS KAKENHI (grant number 19K03608 and 20H00001) and JST CREST (JPMJCR17N2).


\begin{thebibliography}{999}
\bibitem{K-M} A. Kock and J.Mikkelsen, Topos-theoretic factorization of nonstandard extensions, in \textit{Victoria symposium on nonstandard analysis}, 
\textit{Lecture Notes in Mathematics} \textbf{369} (Springer-Verlag, 1974), 122-143.  
\bibitem{ETCS} W. Lawvere, Elementary Theory of the Category of Sets, Proc. Nat. Acad. Sci. \textbf{52} (1964), 1506-1511.
\bibitem{Mac} S. MacLane, \textit{Categories for the working mathematician}, (Springer: Berlin, Heidelberg, New York, 2nd ed. 1998).
\bibitem{Nel} E. Nelson, Internal set theory: A new approach to nonstandard analysis, Bull. Amer. Math. Soc. \textbf{83} (1977), 1165-1198.
\bibitem{Rob} A. Robinson, \textit{Non-standard Analysis}, (North Holland, 1966). 
\bibitem{Roe} J. Roe, \textit{Lectures on Coarse Geometry}, \textit{University Lecture Series} \textbf{31} (American Mathematical Society, Providence, RI, 2003). 

\end{thebibliography}
\end{document}